\documentclass[12pt,thmsa]{article}
\usepackage{amsfonts}
\usepackage[dvips]{graphics}
\usepackage{amsmath}
\usepackage{amssymb}
\usepackage{latexsym}

\def\mylabel#1{\label{#1}}
\newtheorem{theorem}{Theorem}[section]
\newtheorem{lemma}[theorem]{Lemma}
\newtheorem{corollary}[theorem]{Corollary}
\newtheorem{proposition}[theorem]{Proposition}
\newtheorem{example}[theorem]{Example}
\newtheorem{remark}[theorem]{Remark}

\newtheorem{hypothesis}[theorem]{Hypothesis}

\newtheorem{definition}[theorem]{Definition}
\def\bit{\begin{itemize}}
\def\eit{\end{itemize}}
\reversemarginpar   
\def\bc{\begin{center}}
\def\ec{\end{center}}
\def\bthm{\begin{theorem}}
\def\ethm{\end{theorem}}
\def\bcor{\begin{corollary}}
\def\ecor{\end{corollary}}
\def\bprop{\begin{proposition}}
\def\eprop{\end{proposition}}
\def\blem{\begin{lemma}}
\def\elem{\end{lemma}}
\def\bex{\begin{example} {\rm }
\def\eex{\end{example} }}
\def\brem{\begin{remark}}
\def\erem{\end{remark}}
\def\prf{\noindent{\bf Proof. }}
\def\bdes{\begin{description}}
\def\edes{\end{description}}
\def\ita{\item[(a)]}
\def\itb{\item[(b)]}

\def\iti{\item[(i)]}
\def\itii{\item[(ii)]}
\def\itiii{\item[(iii)]}

\def\beq{\begin{equation}}
\def\eeq{\end{equation}}
\def\ben{\begin{enumerate}}
\def\een{\end{enumerate}}
\def\beqar{\begin{eqnarray}}
\def\eeqar{\end{eqnarray}}
\def\beqarr{\begin{eqnarray*}}
\def\eeqarr{\end{eqnarray*}}

\def\RR{{\mathbb R}}  

 \def\cE{\mathcal{E}} \def\cF{\mathcal{F}}
\def\cG{\mathcal{G}}  
  
 \def\cN{\mathcal{N}} 
  \def\cR{\mathcal{R}}

\def\R#1{{\mathbb R}^{#1}}   
\def\Rp{{\mathbb R}_+}   
\def\qed{\hspace{.1in}{$\blacksquare$}}

\def\Pr{{\mathbb P}} 

\def\E{{\mathbb E}} 

\def\NN{{\mathbb N}}       
\def\one{{\bf 1}}
\def\rar{\rightarrow}

\def\eps{\epsilon}
\def\bar(\|)

\def\la{\langle}

\def\dist{\hbox{dist}}

\def\d#1dt{\frac{d#1}{dt}}    

\def\R{\mathbb R}

\def\hop{{\noindent}}

\newcommand{\vs}[1]{\vspace{#1}}

\newcommand{\begitem}{\begin{itemize}}
\newcommand{\finit}{\end{itemize}}

\newcommand{\zun}{\vs{0.1cm}} 

\newcommand{\zdeux}{\vs{0.2cm}} 
\newcommand{\ztrois}{\vs{0.3cm}}
\newcommand{\zcinq}{\vs{0.5cm}} 

\newtheorem{theoreme}{Theorem}[section]

\newtheorem{hypothese}[theoreme]{Hypothesis}

\begin{document}
\title{Stochastic Approximation, Cooperative Dynamics and Supermodular Games\thanks{We acknowledge financial support from the  Swiss National Science Foundation Grant 200021-103625/1}}
\author{
Michel Bena\"{\i}m
\\{\small michel.benaim@unine.ch}
\\ Mathieu Faure
\\ {\small mathieu.faure@unine.ch}
\zdeux
\\{\small Institut de Math\'ematiques, Universit\'e de Neuch\^atel},
\\{\small Rue Emile-Argand 11. Neuch\^atel. Switzerland}.
}

\maketitle

{\small
\begin{abstract} This paper considers a stochastic approximation algorithm, with decreasing step size and martingale difference noise. Under very mild assumptions, we prove the non convergence of this process toward a certain class of repulsive sets for the associated  ordinary differential equation (ODE). We then use this result to derive the convergence of the process when the ODE is {\itshape cooperative} in the sense of \cite{Hir85}. In particular,  this allows us to extend significantly the main result of \cite{HofSan02} on the convergence of {\itshape stochastic fictitious play} in {\itshape supermodular games}.
\end{abstract}
}

\ztrois

\hop {\bfseries MSC2010 Subject classification:} 62L20, 37C50, 37C65, 37Dxx, 91A12.

\section{Introduction}
\mylabel{sec:intro}

Let $F: \mathbb{R}^d \rightarrow \mathbb{R}^d$ be a smooth vector field and $\left(\Omega, \mathcal{F}, \mathbb{P} \right)$ a probability space. We consider a $\mathbb{R}^d$-valued discrete time stochastic process $(x_n)_n$  whose general form can be written as the following recursive formula:

\begin{equation}
\label{eq:DSA}
x_{n+1} - x_n = \frac{1}{n+1} \left( F(x_n) + U_{n+1}\right),
\end{equation}
We assume that $\left(\Omega, \mathcal{F}, \mathbb{P} \right)$ admits a filtration $(\mathcal{F}_n)_n$  such that $x_0$ is measurable with respect to $\mathcal{F}_0$, and $(U_n)_n$ is a $(\mathcal{F}_n)_n$-adapted  sequence of random shocks (or perturbations). Throughout the paper, we make the following assumptions:

\begin{hypothese} \label{RM} We assume that:

\begitem
\item[$(i)$] $(U_n)_n$ is a martingale difference: for any $n \in \mathbb{N}^*$,
\[\mathbb{E} \left(U_{n+1} \mid \mathcal{F}_n\right) = 0.\]
\item[$(ii)$] $F$ is Lipschitz continuous, with Lipschitz constant $L$.
\end{itemize}
\end{hypothese}
\hop Such a stochastic approximation process is generally referred to as a Robbins-Monro algorithm (see \cite{RobMon51} or  \cite{KieWol52}). A natural approach  to obtain information on the asymptotic behavior of the sample paths $(x_n(\omega))_n$ is to compare them to the trajectories of the ordinary differential equation
\begin{equation} \label{ODE}
\dot{x} = F(x).
\end{equation}

\hop Indeed, one can interpret (\ref{eq:DSA}) as some kind of Cauchy-Euler approximation scheme for solving this ODE numerically, with a decreasing step size and an added noise. Since we assume that the noise has null expectation conditionally to the past, it is natural to expect that, for almost every $\omega \in \Omega$, the limit sets of the sample paths $(x_n(\omega))_n$ are related to the asymptotic behavior of the ODE solution curves. This approach was first introduced in \cite{Lju77}  and is usually referred to as the ODE method. Thereafter, the method has been studied and developed by many authors (including \cite{KusCla78},  \cite{BMP90},  \cite{Duf96} or  \cite{KusYin03}) for very simple dynamics (e.g. linear or gradient-like).
\zdeux

\hop In a series of papers (\cite{BenHir96} and \cite{Ben96} essentially), Bena\"{\i}m and Hirsch proved that the asymptotic behavior of $(x_n)_n$ can be described with a great deal of generality  through the study of the asymptotics of (\ref{ODE}), regardless of the nature of $F$. In particular, under certain assumptions on the noise,
\begin{itemize}
\item[$(a)$] the limit sets of $(x_n)_n$ are almost surely \emph{internally chain recurrent} in the sense of Bowen and Conley (see  \cite{Bow75} and \cite{Con78}). This result is detailled in section 2.1.
\item[$(b)$] the random process $(x_n)_n$ converges with positive probability to any given attractor of (\ref{ODE}). See theorem 7.3 in \cite{Ben99} for a precise statement.
\end{itemize}

\hop In addition, it was proved in \cite{Pem90}  that, with probability one,  $(x_n)_n$ does not converge  to linearly unstable equilibria. Some additional references to non convergence results are given in section 3.

The motivation of this paper is threefold. First, under some additional assumptions on the noise, we prove the non convergence of $(x_n)_n$ toward a certain class of unstable sets (including linearly unstable equilibria, periodic orbits and normally hyperbolic  sets), under less regularity assumptions than the existing results. This is detailled in  section 3.

\hop Secondly, in section 4, we use these results, combined with with the nature of limit sets (see point $(a)$ above) and the structure of chain recurrent sets for cooperative dynamics (see  \cite{Hir99}) to prove convergence of $(x_n)_n$ to the set of "stable" equilibria when $F$ is cooperative and irreducible. This  answers a question raised in \cite{Ben00}.

\hop Finally, these results are applied to prove the convergence of {\itshape stochastic fictitious play} in {\itshape supermodular games} in full generality. This proves a conjecture raised in  \cite{HofSan02}.

\section{Background, Notation  and Hypotheses}
\mylabel{sec:background}

Let $F$ denote a locally Lipschitz  vector field on $\RR^d$. By standard results, the Cauchy problem
$\d{y}dt = F(y)$ with initial condition $y(0) = x$ admits a unique solution
$t \rar \Phi_t(x)$ defined on an open interval $J_x \subset \RR$ containing the origin. For simplicity in the statement of our results we furthermore assume that $F$ is {\em globally integrable}, meaning that $J_x = \RR$ for all $x \in \RR^d.$ This holds in particular if $F$ is sublinear; that is
$$\limsup_{||x|| \rar \infty} \frac{||F(x)||}{||x||} < \infty.$$
We let $\Phi = \{\Phi_t\}_{t \in \RR}$ denote the flow induced by $F.$

A continuous map $\chi : \mathbb{R}_+ \rightarrow \mathbb{R}^d$ is called an {\itshape asymptotic pseudo trajectory} (APT) for $\Phi$ \cite{BenHir96}  if, for any $T>0$,
\[\lim_{t \rightarrow + \infty} d_\chi(t,T) = 0,\]
where
\beq
\label{eq:defdX}
d_\chi(t,T) =  \sup_{h \in [0,T]} \left\|\chi(t+h) - \Phi_h(\chi(t)) \right\|.
\eeq
In other terms, for any $T>0$, the curve joining $\chi(t)$ to $\chi(t+T)$ shadows the trajectory of the semiflow starting from $\chi(t)$ with arbitrary accuracy, provided $t$ is large enough.
\brem
\mylabel{T=1}
Assume that $\Phi_1$ restricted to $\chi(\RR_+)$ is uniformly continuous. This holds in particular if $\chi$ or $F$ are bounded maps. Then
$$\lim_{t \rar \infty} d_\chi(t,1) = 0 \Leftrightarrow \forall~T > 0, \lim_{t \rar \infty} d_\chi(t,T) = 0.$$
\erem
Let $(\Omega, \cF, \Pr)$ be a probability space  equipped with   some non decreasing sequence of $\sigma$-algebras $(\mathcal{F}_t)_{t \geq 0}.$ Throughout this paper  we will consider an $(\cF_t)_t$-adapted continuous time stochastic process $X = (X(t))_{t \geq 0}$ verifying the following condition:
\begin{hypothese}\mylabel{H1} {\rm   There exists a map $\omega: \mathbb{R}_+^3 \rightarrow \mathbb{R}_+$ such that:
\bdes
\iti  For any $\delta > 0, T > 0$,
\[\mathbb{P} \left(\sup_{s \geq t} d_X(s,T) \geq \delta \mid \mathcal{F}_t \right)  \leq \omega (t,\delta,T), \]
\itii
$\lim_{t \rar \infty} \omega(t,\delta,T) = 0.$
\edes
}
\end{hypothese}
A sufficient condition  ensuring hypothesis \ref{H1} is that
\beq
\label{eq:H1+}
\mathbb{P} \left( d_X(t,T) \geq \delta \mid \mathcal{F}_t \right)  \leq \int_{t}^{t+T} r(s,\delta,T) ds
\eeq
for some $r : \R^3 \mapsto \Rp$
such that
$$\int_0^{\infty} r(s,\delta,T)ds < \infty.$$
In this case
$$\omega(t,\delta,T) = \int_{t}^{\infty} r(s,\delta,T) ds.$$
The proof of the following proposition is obvious.
\bprop Under  hypothesis \ref{H1}, $X$ is almost surely an asymptotic trajectory for $\Phi.$
\eprop
\bex{\bf (Diffusion processes)}
\mylabel{ex:diffusion}
{\rm

Let $X$ be solution to the stochastic differential equation
\[dX(t) = F(X(t)) dt + \sqrt{\gamma(t)} dB_t,\]
where $F$ is a globally Lipschitz vector field, $(B_t)$ a standard Brownian motion on $\RR^d$ and $\gamma : \Rp \mapsto \Rp$ a decreasing continuous function.
Assume that
\[\int_0^{+ \infty} \exp \left(\frac{-c}{\gamma(t)} \right) dt < + \infty\]
for all $c > 0.$ Then (\ref{eq:H1+})
 is satisfied with
\[r(t,\delta,T) =  C \exp \left( - \frac{\delta^2 C(T)}{\gamma(t)} \right )\]
 where $C$ and $C(T)$ are positive constants. This is proved in (\cite{Ben99}, Proposition 7.4)
 }
\eex

\bex{\bf (Robbins-Monro algorithms)}
\mylabel{ex:RM}{\rm Let $(x_n)_n$ be a stochastic approximation algorithm governed by the recursive formula

\begin{equation}
x_{n+1} - x_n = \gamma_{n+1} \left( F(x_n) + U_{n+1}\right),
\end{equation}
where $\gamma_n \geq 0, \, \sum_n \gamma_n = \infty,$
and which satisfies Hypothesis \ref{RM}. Assume  furthermore that
 one of the two following conditions holds:
\begitem
\iti  There exists some $q \geq 2$ such that
\[\sum \gamma_n^{1+ q/2} < + \infty \, \mbox{ and } \; \, \sup_n \mathbb{E}\left(\|U_n\|^q \right) < + \infty.\]
\itii \bdes
\ita The sequence $(U_n)_n$ is {\em subgaussian} (for instance bounded) meaning that
$$\E(\exp (\langle \theta, U_{n+1} \rangle)) | \cF_n) \leq \exp (\Gamma ||\theta||^2)$$ for some $\Gamma > 0;$ and
\itb  for any $c> 0$,
\[\sum_{n} \exp \left(\frac{-c}{\gamma_n} \right) < + \infty.\]
\edes
\end{itemize}

Set $\tau_n := \sum_{i=1}^n \gamma_i$. We call $X$ the continuous time affine interpolated process induced by $(x_n)_n$ and $\overline{\gamma}$ the piecewise constant deterministic process induced by $(\gamma_n)_n$:
\[X (\tau_i + s) := x_i + s \frac{x_{i+1} - x_i}{\gamma_{i+1}}, \mbox{ for } \; i \in \mathbb{N}, \; \, s \in [0,\gamma_{i+1}]\]
  and
\[\overline{\gamma}(\tau_i + s) := \gamma_{i+1} \mbox{ for } \; s \in [0, \gamma_{i+1}[.\]
Under one of the above condition $(i)$ or $(ii)$, this continuous time process is an asymptotic pseudo trajectory of the flow induced by $F$ (see \cite{Ben99}). Additionally,  we have the following result  (see \cite{Ben99}  and more specifically \cite{Ben00}):

\begin{proposition} \label{omega}
Let $k_0 := \inf \left\{ k \in \mathbb{N} \mid \; \, \gamma_k \leq \frac{B \delta^2}{2} \right\}$. Then, for any $s \geq \tau_{k_0}$, condition (\ref{eq:H1+}) holds with
$$r(s,\delta,T) =  \frac{B  \overline{\gamma}^{q/2}(s)}{\delta^q}  $$
in the first case, and
$$r(s,\delta,T) =  2d  \exp \left( \frac{-B \delta^2}{\overline{\gamma}(s)}\right)$$
in the second, where $B$  is some positive constant depending only on the noise, the step size and the vector field.

\end{proposition}

}
\eex

\subsection{The Limit set Theorem}
A set $L \subset \RR^d$ is said to be {\em invariant} (respectively {\em positively invariant}) for $\Phi$ provided $\Phi_t(L) \subset L$ for all $t \in \RR$ (respectively $t \in \Rp$).

Let $L$ be an  invariant set for $\Phi.$ We let $\Phi^L$ denote the restriction of $\Phi$ to $L.$ That is,
$\Phi_t^L(x) = \Phi_t(x)$ for all $x \in L$ and $t \in \RR.$
Note that with such a notation $\Phi = \Phi^{\RR^d}.$

An {\em attractor } for $\Phi^L$ is a nonempty compact invariant set $A \subset L$ having a neighborhood $U$ in $L$ such that
$$\lim_{t \rar \infty} \dist(\Phi^L_t(x), A) = 0$$ uniformly in $x \in U.$
Note that if $L$ is compact, $L$ is always an attractor for $\Phi^L.$ An attractor for $\Phi^L$ distinct from $L$ is called a {\em proper attractor}.

The {\em basin of attraction} of $A$ for $\Phi^L$ is the open set (in $L$) consisting of every $x \in L$ for which $\lim_{t \rar \infty} \dist(\Phi_t(x),A) = 0.$

A {\em global attractor} for $\Phi$ is an attractor which basin is $\RR^d.$
If such an attractor exists, $\Phi$ (respectively $F$)  is called a {\em dissipative} flow (respectively vector field).

A compact invariant set $L$ is said to be {\em internally chain-transitive} or {\em attractor free} if $\Phi^L$ has no proper attractor (see e.g. \cite{Con78}).

A fundamental property of asymptotic pseudo trajectories is given by the following result due to \cite{Ben96}  for stochastic approximation processes and \cite{BenHir96}  for APT. We refer to  \cite{Ben99} for a proof and more details; and also to \cite{Pem07} for a recent overview and some applications.
\bthm
\mylabel{th:BH} Let $\chi$ be a bounded APT, then its limit set $$\mathcal{L}(\chi) = \bigcap_{t \geq 0} \overline{\chi([t,\infty[)}$$ is internally chain transitive.
\ethm
\bcor
\mylabel{th:BH2} Under hypothesis \ref{H1}, the limit set of $X$ is
almost surely internally chain transitive on the event $\{\sup_{t \geq 0 } \|X(t)\| < \infty\}.$
\ecor
\section{Non convergence toward normally hyperbolic repulsive sets}
\mylabel{sec:nonconv}
From Corollary \ref{th:BH2} we know that the limit set of $X$ is internally chain transitive (ICT). However  not every ICT set can  be such a limit set because the noise may  push the process away from certain ``unstable'' sets.
For  equilibria this  question has been tackled by several authors including \cite{Pem90}, \cite{Tar}, \cite{BraDuf96} and it was proved that, under natural conditions, $X$ has zero probability to converge toward a linearly unstable equilibrium.  This  has been extended to linearly unstable periodic orbit by \cite{BenHir95} and to more general normally hyperbolic sets by \cite{Ben99}). The proofs of all these results rely on the assumption that  the unstable manifold of the set (to be defined below) is sufficiently smooth (at least $C^{1+\alpha}$ with $\alpha > 1/2).$ While for linearly unstable equilibria or periodic orbit such a  regularity  assumption  follows directly  from the regularity of the vector field,  the situation is much trickier for more general sets.

The purpose of this section is to extend the non convergence results mentioned above under less regularity assumptions. This will prove to be of fundamental importance in our analysis of cooperative dynamics and supermodular games in  section  \ref{sec:coop}.

Let $S$ be a $\mathcal{C}^1$, $(d-k)$-dimensional ($k \in \{1,..,d\}$) submanifold of $\RR^d$ and $\Gamma$  a compact  invariant set contained in $S$. We assume that $S$ is {\em locally invariant} meaning that there exists a neighborhood $U$ of $\Gamma$ in $\RR^d$  and a positive time $t_0$ such that
$$\Phi_t(U \cap S) \subset S$$ for all $|t| \leq t_0.$
We let $\mathcal{G}(k,d)$ denote  the Grassman manifold of $k$ dimensional planes in $\RR^d$.  For $p \in S$, the tangent space of $S$ in $p$ is denoted $T_p S$.

\begin{definition}
{\rm  $\Gamma$ is  called a {\em normally hyperbolic repulsive set} if there exists a continuous map

\[ p \in \Gamma \mapsto E^u_p \in \mathcal{G}(k,d),\]

such that
\begitem

\iti for any $p \in \Gamma$,

\[\mathbb{R}^d = T_p S \oplus E^u_p,\]

\itii for any $t \in \mathbb{R}$ and any $p \in \Gamma$,

\[D \Phi_t(p) E^u_p = E^u_{\Phi_t(p)},\]

\itiii there exists positive constants $\lambda$ and $C$ such that, for any $p \in \Gamma$, $w \in E^u_p$ and $t \geq 0$, we have

\[\left\|D \Phi_t(p) w \right\| \geq C e^{\lambda t} \|w\|. \]

\end{itemize}

}
\end{definition}
The two basic examples of normally hyperbolic sets are the following. For more details, see (\cite{Ben99}, Section 9).

\bex {\rm{\bfseries (Linearly unstable equilibrium):} If $\Gamma = \{p\}$, where $p$ is a linearly unstable equilibrium (not necessarily hyperbolic), then it is a normally hyperbolic repulsive set.}
\eex

\bex {\rm {\bfseries (Hyperbolic linearly unstable periodic orbit):} If $\Gamma$ is a periodic orbit, the unity is always a Floquet multiplier. It is hyperbolic if the others multipliers all have moduli different from $1$ and it is linearly unstable if at least one has modulus strictly greater than one. If both assumptions are checked then $\Gamma$ is a normally hyperbolic repulsive set.}
\eex
For further analysis, it is convenient to extend
 the map $p \rar E^u_p$ to a neighborhood of $\Gamma$ and to  approximate it  by a smooth map. More precisely it is shown in \cite{Ben99}, Section 9.1 that
there exists
 a neighborhood $\cN_0 \subset U$ of $\Gamma$ and a $C^1$ bundle $$\tilde{E}^u = \{(p,v) \in S \cap \cN_0 \times \RR^d:  \: v \in \tilde{E^u_p}\}$$
where $\tilde{E}^u_p \in \cG(k,d)$ such that:
\bdes
\iti For all $p \in S \cap \cN_0$, $\RR^d = T_p S \oplus \tilde{E}^u_p;$
\itii the map $H : \tilde{E}^u \mapsto \RR^d$ defined by $H(p,v) = p+v$ induces  a $C^1$ diffeormorphism from a neighborhood of the zero section $\{(p,0) \in \tilde{E}^u\}$ onto $\cN_0.$
\edes
Let now  $V : \cN_0 \mapsto \R_+$ be the map defined by
 $V(x) = ||v||$ for $H^{-1}(x) = (p,v).$
 The form of $V$  implies  that there exits $L > 0$ such that
\beq
\label{Vlip}
d(x,S) \leq V(x) \leq L d(x,S)
\eeq
 for all $x \in \cN_0.$ Then according to Lemma 9.3 in \cite{Ben99}
there exist a bounded  neighborhood $\cN_1 \subset \cN_0$ of $\Gamma,$
 and numbers $T > 0$, $\rho > 1$ such that
\beq
\label{Vexp}
\forall x \in \cN_1, \; \,  V(\Phi_T(x)) \geq \rho V(x).
\eeq

Given a  neighborhood $\cN \subset U$ of $\Gamma$ we let
\[\mathsf{Out}_{\eps} = \mathsf{Out}_{\eps}(\cN,S) := \left\{ x \in \cN  \mid d(x, S \cap\cN) \geq \eps \right\}.\]
and
\[\mathsf{In}_{\eps} = \mathsf{In}_{\eps}(\cN,S) := \cN \setminus \mathsf{Out}_{\eps}.\]
\blem
\label{th:out}
\bdes

\iti
There exists a bounded neighborhood $\cN \subset U$ of $\Gamma, T > 0$ and $\rho > 1$ such that for all $\eps > 0,$
$$\Phi_T(\mathsf{Out}_{\eps}(\cN,S)) \cap \cN \subset \mathsf{Out}_{\rho \eps}(\cN,S).$$
In particular, every compact invariant subset contained in $\cN$ lies in $S.$
\itii For all $R > 0$ there exists a finite set $\{v_1,\ldots,v_n\} \subset \R^d$ and a Borel map $I : \Gamma \mapsto \{1,\ldots,n\}$ such that for all $p \in \Gamma$ and $v \in B(v_{I(p)}, 1),$
$$p + \eps v \in \mathsf{Out}_{R\eps}.$$
\edes
\elem
\prf
Choose $k \in \NN$ such that $\rho^k > L$ and $\cN \subset \cN_1$ be small enough so that $\Phi_{kT}(\cN) \subset \cN_1.$ Then, using (\ref{Vlip}) and (\ref{Vexp})
for all $x \in \cN$,
$$d(\Phi_{kT}(x),S) \geq \frac{1}{L} V(\Phi_{kT}(x)) \geq \frac{\rho^k}{L} V(x)
\geq  \frac{\rho^k}{L} d(x,S).$$
Replacing $T$ by $kT$ and $\rho$ by $\frac{\rho^k}{L}$ gives the result.

We now prove the second assertion. Given $R > 0,$ let $f : \Gamma \mapsto \RR^d$ be a measurable function such that for all $p \in \Gamma,  f(p) \in \tilde{E}^u_p$ and $\|f(p)\|= L(R+2)$ where $L$ is the constant appearing in  (\ref{Vlip}).  The bundle $\tilde{E}^u$ being locally trivial, it is not hard to construct such a function. By compactness of $\overline{f(\Gamma)}$, there exists a finite set $\{v_1,\ldots,v_n\} \subset f(\Gamma)$ such that $f(\Gamma)\subset \cup_{i = 1}^n B(v_i,1).$
For $p \in \Gamma,$ set  $$I(p) = \min \{i = 1, \ldots, n\, : \|f(p) -v_i\|\leq 1\}.$$
Then, for $I(p) = i$ and $v \in B(v_i,1)$,
$$d(p + \eps f(p),S) \leq d(p+ \eps v,S) + \eps \|f(p) -v)\| \leq d(p+\eps v,S) + 2 \eps.$$
On the other hand, by (\ref{Vlip}), $$d(p + \eps f(p),S) \geq \frac{1}{L} V(p+ \eps f(p)) = \frac{\eps \|f(p)\|}{L} = \eps(R+2).$$
Hence
$$d(p+\eps v,S) \geq R \eps.$$
\qed

\bcor \label{th:nonconv0}
Let $\cN, T$ and $\rho$ be like in Lemma \ref{th:out}, and set $\delta = (\rho -1) > 0.$ Let $Y$ be an asymptotic pseudo-trajectory verifying
\bdes
\iti $\chi(0) \in \mathsf{Out}_{\eps},$
\itii for all $t \geq 0, d_\chi(t,T) \leq \delta \eps.$
\edes
Then $\chi$ eventually leaves $\cN.$
\ecor
\prf Suppose that $\chi$ remains in $\cN.$ We claim that $\chi(kT) \in \mathsf{Out}_{\eps}$ for all $k \in \NN.$  If
 $\chi(kT) \in \mathsf{Out}_{\eps}$ then $\Phi_T(\chi(kT)) \in \mathsf{Out}_{\rho \eps}$ by Lemma \ref{th:out}. Hence $\chi(kT+T) \in \mathsf{Out}_{\eps}$ since $d_\chi(kT,T) \leq \delta \eps.$ This proves the claim by induction on $k.$
It follows that the limit set of $\chi$ meets $\mathsf{In}_{\eps}$ but, by the limit set theorem \ref{th:BH} and Lemma \ref{th:out}, this limit set has to be in $S.$ A contradiction.
\qed

\subsection{Non convergence: sufficient conditions}
Throughout this section we let $\cN, T$ and $\rho$ be like in Lemma \ref{th:out}, and  $\delta = (\rho -1) > 0.$
We  let $X$ be a continuous time $(\cF_t)$-adapted process verifying hypothesis \ref{H1} and $E_t$ be the event
\[E_t = \{\forall s \geq t: \:   X(s) \in \cN\}.\]
\begin{lemma}
\label{th:nonconv1} On the event $\{X(t) \in \mathsf{Out}_{\eps}\}$,
$$\Pr(E_t | \cF_t) \leq \omega(t, \delta \eps, T)$$
 and
 $$\Pr(E_t  | \cF_t) \leq 1 - [1- \omega(t+1,\delta \eps,T)]
 \Pr(X(t+1) \in \mathsf{Out}_{\eps}  |  \cF_t).$$
\end{lemma}
\prf
The first inequality   follows from Corollary \ref{th:nonconv0}.
Now
\begin{eqnarray*}
&& \Pr(E_t | \cF_t) \leq \Pr(E_{t+1} | \cF_t) \\
&=& \Pr(E_{t+1}; X(t+1) \in \mathsf{Out}_{\eps}|\cF_t) + \Pr(E_{t+1}; X(t+1) \in \mathsf{In}_{\eps}|\cF_t)\\
&=& \E( \Pr(E_{t+1}|\cF_{t+1}) \one_{X(t+1) \in \mathsf{Out}_{\eps}}|\cF_t) + \E( \Pr(E_{t+1}|\cF_{t+1}) \one_{X(t+1) \in \mathsf{In}_{\eps}}|\cF_t)\\
&\leq& \omega(t+1,\delta \eps, T) \Pr (X(t+1) \in \mathsf{Out}_{\eps}|\cF_t) + \Pr (X(t+1) \in \mathsf{In}_{\eps}|\cF_t).
\end{eqnarray*}
 \qed
\begin{lemma} \label{critere}
Assume that there exists a maps $\eps : \R_+ \mapsto \R_+$
with $\lim_{t \rar \infty} \eps(t) = 0$ and constants $c > 0$ and $c' < 1$ such that for $t$ large enough
\bdes
\iti $\Pr (X(t+1) \in \mathsf{Out}_{\eps(t)}|\cF_t) \geq c$ on the event  $\{{X(t) \in \mathsf{In}_{\eps(t)}}\}.$
\itii  \[  \omega \left(t, \delta \eps(t), T \right)  < c', \]

\edes
Then
$$\Pr(X(t) \rar \Gamma) = 0.$$
\end{lemma}
\prf
One has $$\{X(t) \rar \Gamma\} \subset \bigcup_{n \in \NN} E_n$$ and it suffices  to prove that $\Pr(E_n) = 0$ for all $n \in \NN.$

For all $t \geq n$, $E_n \subset E_{t}.$ Thus  $$\Pr(E_n | \cF_{t}) \leq \Pr(E_{t}|\cF_{t}) \leq \max(c', 1-(1-c')c ),$$
where the last inequality follows from the assumptions and   Lemma \ref{th:nonconv1}.
Now, by a classical Martingale result, $$1 > \max(c', 1-(1-c')c ) \geq \lim_{t \rar \infty} \Pr(E_n  |\cF_{t}) \rar \one_{E_{n}}$$ almost surely.
Hence the result.
\qed

\begin{hypothesis}
\mylabel{H2}  Assume that there exists a map $\gamma : \Rp \rar \Rp$ with $\lim_{t \rar \infty} \gamma(t) = 0$ and an adapted process $(Y(t))_{t \geq 0}$ such that
\bdes
\iti For all $\eps > 0,$ $$\lim_{t \rar \infty} \Pr \left( \left\| \frac{X(t+1) - \Phi_1(X(t))}{\sqrt{\gamma(t)}} - Y(t+1)\right\| \geq \eps | \cF_t \right) = 0,$$
\itii For all open set $O \subset \R^d$
$$ \liminf_{t \rar \infty} \Pr( Y(t+1) \in O | \cF_t) > 0.$$
\itiii There exists $a > 0$ such that  $$\limsup_{t \rar \infty} \omega(t, a \sqrt{\gamma(t}),T) < 1.$$
\edes
\end{hypothesis}
\bthm
\mylabel{th:nonconvX} Let $X$ be a continuous $(\cF_t)$-adapted process verifying hypotheses \ref{H1} and \ref{H2}. Then $$\Pr(X(t) \rar \Gamma) = 0.$$ \ethm
\prf We shall prove that the assumptions of Lemma \ref{critere} are fulfilled with  $\eps(t) = \frac{\sqrt{\gamma(t)}}{\alpha};$ where $\alpha = \frac{\delta}{a}$ and $a$ is given by hypothesis \ref{H2}(iii).
Condition $(ii)$ of the lemma is clearly verified.

To check condition $(i)$
we assume  that $X(t) \in \mathsf{In}_{\eps(t)}$. Hence (for $t$ large enough), $\Phi_1(X(t))$ lies in  $\cN_0 \subset \cN$ and we can
write $$\Phi_1(X(t)) = p(t) + v(t)$$ with $(p(t),v(t)) \in \tilde{E}^u_{p(t)}$ (see the beginning of the section).
Then, by the triangle inequality,
\begin{eqnarray*}
d(X(t+1), S) &\geq&  d(p(t) + \alpha \eps(t) Y(t+1),S) - \|v(t)\| \\
&& - \alpha \eps(t) \|Y(t+1) - \tilde{Y}(t+1)\|.\end{eqnarray*}
with $$\tilde{Y}(t+1) = \frac{ X(t+1) - \Phi_1(X(t)) } {\alpha \eps(t)}.$$
 Now $$\|v(t)\| = V(\Phi_1(X(t))) \leq L d(\Phi_1(X(t),S)) \leq M \eps(t)$$
where the first inequality follows from the Lipschitz continuity of the map $V$ (see (\ref{Vlip})), and the second from the Lipschitz continuity of $\Phi_1$ and invariance of $S$. Thus
$$\frac{d(X(t+1),S)}{\eps(t)} \geq U_t - V_t - M$$
where
$$U_t = \frac{d(p(t) + \alpha \eps(t) Y(t+1),S)}{\eps(t)}$$ and $$V_t = \alpha \|Y(t+1) - \tilde{Y}(t+1)\|.$$

Let $R = \frac{2 + M}{\alpha}.$
Then by lemma \ref{th:out} (ii) and hypothesis \ref{H2} (ii), there exists $c > 0$ such that
$$\Pr (U_t \geq (1 + M)  | \cF_t) =  \Pr(p(t) + \alpha \eps(t) Y(t+1) \in \mathsf{Out}_{R \alpha \eps(t)}|\cF_t) \geq 2c.$$
Furthermore, by Hypothesis \ref{H2},
$$\lim_{t \rar \infty} \Pr (V_t \geq 1 | \cF_t) \leq c$$ for $t$ large enough.
It follows that
$$\Pr\left(\frac{d(X(t+1), S)}{\eps(t)} \geq 1 |\cF_t \right) \geq \Pr (U_t - V_t \geq M+1|\cF_t)$$
$$ \geq \Pr(U_t \geq 2 + M|\cF_t) - \Pr(V_t \geq 1 |\cF_t) \geq c.$$ This proves that condition $(i)$ of the lemma is verified.
\qed

\begin{proposition}
\mylabel{th:nonconvsde} Let $X$ be like in example \ref{ex:diffusion}. Set $l(t) = \log(\gamma(t)).$ Assume that
\bdes
\iti Function $l$ is sub-additive: $l(t+s) \leq l(t) + l(s)$. This holds in particular if $l$ is concave and $l(0) = 0$
\itii There exist constants $a \geq b > 0$ such that
 $- a \leq \dot{l}(t) \leq -b.$
 \edes
 Then hypothesis \ref{H2} holds. In particular, conclusions of Theorems \ref{th:BH} and \ref{th:nonconvX} hold. \end{proposition}
The proof is given in appendix.

We now apply these results to the specific case of Robbins-Monro algorithm. An additional assumption on the noise is needed:

\begin{hypothesis} \label{noise}
There exists positive real values $0 < \Lambda^- < \Lambda^+ < + \infty$ and a continuous map
\[Q: \mathbb{R}^d \rightarrow \mathcal{S}^+\left(\mathbb{R}^d\right) \cap [\Lambda^- I_d, \Lambda^+ I_d],\]

\hop such that $\mathbb{E} \left(U_{n+1} U_{n+1}^T \mid \mathcal{F}_n \right) = Q(x_n)$.

\end{hypothesis}

\begin{proposition} \label{RMNC} Let $(x_n)_n$ be a Robbins-Monro algorithm like in example \ref{ex:RM}  with $\gamma_n = 1/n$ and $\mathbb{E} \left( \|U_n\|^{2p} \mid \mathcal{F}_{n-1}\right)$ almost surely bounded for some $p>1$,   which noise also satisfies hypothesis \ref{noise}. Then the associated interpolated process $X(t)_{t \geq 0}$ satisfies Hypothesis \ref{H2} and therefore,
\[\mathbb{P} \left(X(t) \rightarrow \Gamma \right) = 0.\]

\end{proposition}
\

\hop {\bfseries  Proof.} In appendix.

\section{Application to cooperative dynamics}
\mylabel{sec:coop}
Throughout this section we assume  that for all $x \in \RR^d$ the Jacobian matrix $DF(x) = (\frac{\partial F_i}{\partial x_j}(x))$ has nonnegative off-diagonal entries and is irreducible. Such a vector field $F$ is said to be {\em cooperative} and {\em irreducible} \cite{Hir85}. We refer the reader to \cite{HirSmi06} for a recent survey on the subject.  We furthermore assume that $F$ is  {\em dissipative}, meaning that it admits a global attractor.

For $x,y \in \mathbb{R}^d$,  $x \geq y$ means that $x_j \geq y_j$ for all $j$. If, additionally, $x \neq y$, we write  $x > y.$ If $x_j > y_j$ for all $j$, it is denoted $x \gg y$. Given two sets $A,B \subset \RR^d$ we write $A \leq B$ provided $x \leq y$ for all $x \in A$ and $y \in B.$ Set $A$ is called {\em unordered}
if for all $x,y \in A, \,  x \leq y \Rightarrow x = y.$

The vector field $F$ being  cooperative and irreducible, its flow has positive derivatives \cite{Hir85}, \cite{HirSmi06}. That is  $ D \Phi_t(x) \gg 0$ for $x \in \RR^d$ and $t > 0.$ This implies that it is {\em strongly monotonic} in the sense that $ \phi_t(x) \gg \phi_t(y)$ for all $x > y$ and $t > 0.$

We let $\cE$ denote the equilibria set of $F.$ A Point  $p \in \cE$  is called {\em linearly unstable} if the Jacobian matrix $DF(p)$ has at least one eigenvalue with positive real part. We let $\cE^+$ denote the set of such equilibria and $\cE^{-} = \cE \setminus \cE^{+}.$

An equilibrium point $p \in \cE$
is said to be {\em asymptotically stable from below} if there exists $x<p$ such that $\phi_t(x) \rightarrow p$. The subset of equilibria which satisfy this property is denoted $\mathcal{E}_{asb}$. Note that if $p \in \mathcal{E}_{asb}$, then there exists a non empty open set of initial conditions from which the solution trajectories converge to $p$. In particular $\mathcal{E}_{asb}$ is countable. Given $p \in \mathcal{E}_{asb}$, we introduce the set of points whose limit set dominates $p$:
\[V(p) := \left\{x \mid \omega (x) \geq p \right\}\]
and we let $S_p$ denotes it  boundary: $S_p := \partial V(p)$.
 The following proposition is basically due to (\cite{Hir88}, Theorem 2.1) but for the $C^1$ regularity proved by \cite{Ter96}. Our statement follows from  Proposition 3.2 in \cite{Ben00}, where more details can be found.

\begin{proposition} There exists a unique equilibrium $p^* \in \mathcal{E}_{asb}$ such that $V(p^*) = \mathbb{R}^d$. For any other $p \in \mathcal{E}_{asb}\setminus
\{p^*\},  S_p$ is a $C^1$  unordered invariant hypersurface diffeomorphic to $\RR^{d-1}.$
\end{proposition}
For $p \in  \mathcal{E}_{asb}\setminus
\{p^*\}$ we let $\cR(\Phi^{S_p})$ denote the chain recurrent set of $\Phi$ restricted to $S_p;$  or equivalently, the union of all internally chain transitive sets contained in $S_p.$ We also set
$$\cR'_p = \cR(\Phi^{S_p}) \setminus \{\cE^{-} \cap S_p\}.$$

The first part of the next Theorem is proved in \cite{Ben00} (see the proof of Theorem 2.1) and the second part restates Theorem 3.3 in
the same paper (relying heavily on \cite{Hir99}.
\bthm
\mylabel{th:coopcr}
For any $p \in \cE_{asb} \setminus \{p^*\}$ the set $\cR'_p$ is a repulsive normally hyperbolic set (in the sense of section \ref{sec:nonconv}).
 Any internally chain transitive set is either an ordered arc included  in $\cE^{-}$  or is contained in  $\cR'_p$ for some $p \in \mathcal{E}_{asb} \setminus \{p^*\}$.
\ethm
\brem By a result of \cite{Jia91}, if $F$ is real analytic, it cannot have a nondegenerate ordered arc of equilibria \erem

As a consequence of these results we get the following
\bthm
\mylabel{th:coopnoise}
Let $X$ be a continuous $(\cF_t)$-adapted stochastic process verifying hypotheses \ref{H1} and \ref{H2}.  Then the limit set of $X$ is almost surely an ordered arc contained in $\cE^{-}.$ In case $F$ is real analytic, $X(t)$ converges almost surely to an equilibrium $p \in \cE^{-}.$
\ethm
\prf Follows from Theorems \ref{th:BH}, \ref{th:coopcr} and \ref{th:nonconvX} \qed
\bcor
 Let $X$ be the process given in example \ref{ex:diffusion} with $- a \leq \frac{\dot{\gamma}(t)}{\gamma(t)} \leq - b$ with $ a \geq  b > 0.$
  Then the conclusions of Theorem \ref{th:coopnoise} hold. \ecor
\bcor
\mylabel{th:cooprm}
Let $(x_n)$ be the Robbins Monro algorithm given in example \ref{ex:RM} with $\gamma_n = \frac{1}{n}$.
Assume that hypothesis \ref{noise} holds. Then  the conclusions of Theorem \ref{th:coopnoise} hold. \ecor

\section{Perturbed best response dynamic in supermodular games}
\zcinq

\subsection{ General settings}
\ztrois

\hop Let us consider a $N$ persons game in normal form. Player i's action set is finite and denoted $A^i$, $\Delta^i$ is the mixed strategies set:
\[\Delta^i := \left\{x^i=(x^i(\alpha))_{\alpha \in A^i} \mid  x^i(\alpha) \geq 0, \; \, \sum_{\alpha \in A^i}  x^i(\alpha) = 1\right\}\]
and  $u^i : A^i \mapsto \RR$ his utility function. The set of action profiles (respectively mixed strategy profiles) is denoted $A := \times_{i=1}^N A^i$ (resp. $\Delta := \times_{i=1}^N \Delta^i$). The utility functions $(u^i)_{i=1,..,N}$ are defined on $A$ but linearly extended to $\Delta$:
\[x=(x^1,..,x^N) \in \Delta \mapsto u^i(x) := \sum_{a=(a^1,..,a^N) \in A} u^i(a) x^1(a^1)...x^N(a^N).   \]
We call $G(N,A,u)$ the game induced by these parameters. Throughout our study, we assume that agents play repeatedly and independently.  By this, we  mean that, denoting $a_n = (a^1_n,..,a^N_n)$ the action profile realized at stage $n$ and $(\mathcal{F}_n)_n$ an adapted filtration, we have
\[\mathbb{P}\left(a_{n+1} = (a^1,...,a^N) \mid \mathcal{F}_n \right) = \prod_{i=1}^N \mathbb{P} \left(a^i_{n+1} = a^i \mid \mathcal{F}_n\right).\]

\hop For $a = (a^1,..,a^N)$, $\delta_{a^i}$ denotes the vertex of $\Delta^i$ corresponding to the pure strategy profile $a^i$ and $\delta_a$ is the extreme point of the polyhedron $\Delta$ relative to the pure strategy profile $a$. At last, $\overline{x}_n$ is the empirical distribution of moves up to time $n$ :
\[\overline{x}_n := \frac{1}{n} \sum_{m=1}^n \delta_{a_m} = \left(\frac{1}{n} \sum_{m=1}^n \delta_{a^1_m},.., \frac{1}{n} \sum_{m=1}^n \delta_{a^N_m}\right).\]

\paragraph{Standing Notation} As usual in game theory we let $a^{-i} = (a^j)_{j \neq i},$
$x^{-i} = (x^j)_{j \neq i}, A^{-i} = \times_{j \neq i}  A^j$ etc.  We may  write $(a^i,a^{-i})$ for $a = (a^1, \ldots, a^N)$ and so on.
\subsection{Perturbed best response dynamic}
To shorten notation let us take the point of view of player $1$.
A {\em choice function} for player $1$ is a continuously differentiable map $C : \RR^{A^1} \mapsto \Delta^1.$

Let $f : \R^{A^1} \mapsto \R^{+}$ be a strictly positive probability density and $\varepsilon \in \R^{A^1}$ a random variable having distribution $f(x)dx.$
We say that $C$ is a {\em good stochastic choice function } if it is induced by such a stochastic perturbation $\varepsilon$, in the following sense: for all $\Pi \in \R^{A^1}$ $C(\Pi)$ is the law of the random variable
$$\mathsf{argmax}_{\beta \in A^1} \left( \Pi(\beta) + \varepsilon(\beta) \right).$$
A classical example of good stochastic choice function is the Logit map:
$$L(\Pi)(\alpha) = \frac {\exp{(\eta^{-1} \Pi(\alpha))}}{\sum_{\beta \in A^1} \exp{(\eta^{-1} \Pi(\beta))}}.$$
It is induced by a stochastic perturbation with extreme value density (see \cite{FudLev98} and  \cite{HofSan02}).

Given a choice function $C$, the {\em smooth} or {\em perturbed best response} associated to  $C$  is the  map $\mathbf{br}^1 : \Delta^{-1} \mapsto \Delta^1$
defined by
 $$\mathbf{br}^1(y) = C(u^1(\cdot~,y)).$$
\begin{definition}
Let $\mathbf{br}^1$ be a perturbed best response for player $1.$ A smooth fictitious play (SFP) strategy induced by $\mathbf{br}^1$ is a strategy such that, for any other opponent's strategy,
\begin{equation}
\label{eq:defbr}
\Pr(a_{n+1}^1 = .  \mid \mathcal{F}_n) = \mathbf{br}^1 (\overline{x}^{-1}_n),
\end{equation}
where  $\overline{x}^{-1}_n$  is the empirical moves of the opponents up to time $n$.
\end{definition}

Stochastic fictitious play was originally  introduced by Fudenberg and Kreps (see \cite{FudKre93}) and the concept behind is that players use fictitious play strategies in a game where payoff functions are perturbed by some random variables in the spirit of \cite{Har73}.
To be more precise, suppose that at time $n+1,$ the payoff function to player $1$ is the map $$u^1_{n+1} : A \mapsto \R,$$
$$a \mapsto  u^1(a) + \varepsilon_{n+1}(a^1),$$
where $\varepsilon_n \in \R^{A^1}$ is a random vector which conditional law, given $\cF_n$ is $f(x)dx.$
Suppose furthermore that $u^1_{n+1}$ is known to player $1$ as well as all the  actions $a_1, \ldots, a_n$ played up to time $n.$
Fictitious play assumes that player $1$ chooses the best response to  $\overline{x}^{-1}_n.$ That is
$$a^1_{n+1} = \mathsf{argmax}_{\beta \in A^1} u^1_{n+1}(\beta,\overline{x}^{-1}_n).$$ Hence equation
(\ref{eq:defbr}) holds where $\mathsf{br}^1$ is the smooth best response associated to the good stochastic choice function induced by $\varepsilon_{n+1}$.\footnote{Another approach is to consider that the player chooses to randomize slightly its moves playing a best response relative to a  payoff function perturbed by a deterministic map. Hofbauer and Sandholm (see \cite{HofSan02}) proved that any admissible stochastic perturbation can be represented in term of a deterministic perturbation.}
On the subject, see also the papers \cite{FudLev95},  \cite{FudLev98} or \cite{BenHir99}.

Let us get back to the settings described earlier with $N$ players. We are interested in the asymptotic behavior of the sequence $\overline{x}_n$ when every player adopts a smooth fictitious play strategy. In the remaining of the section, an $N$-uple of  perturbed best response maps is  given and we let  $\mathbf{br} : \Delta \mapsto \Delta$ denote the map defined by
\[\mathbf{br}(x):= \left(\mathbf{br^1}(x^{-1}),..,\mathbf{br^N}(x^{-N}) \right).\]
The set of {\itshape perturbed Nash equilibria}, i.e. the set of $x \in \Delta$ such that $\mathbf{br}(x) = x$ (which can be viewed as the Nash equilibria in an auxiliary perturbed game) will be refered to as PNE. A simple computation gives
\[\overline{x}_{n+1} - \overline{x}_n = \frac{1}{n+1} \left( \delta_{a_n} - \overline{x}_n \right).\]

\hop Hence, the expected increments satisfy:
\[\mathbb{E} \left(\overline{x}_{n+1} - \overline{x}_n \mid \mathcal{F}_n\right) = \frac{1}{n+1} \left(\mathbf{br}(\overline{x}_n) - \overline{x}_n \right).\]

\hop The recursive formula describing the evolution of the random process $(\overline{x}_n)_n$ can then be written
\begin{equation}
\label{eq:SFP}
\overline{x}_{n+1} = \overline{x}_n + \frac{1}{n+1} \left(F(\overline{x}_n) + U_{n+1} \right),
\end{equation}
where

\begitem

\iti the vector field $F$ defined by $F(x) = \mathbf{br}(x) - x$ is smooth,
\itii the noise $U_{n+1}$ is a bounded martingale difference by construction and given by
\[U_{n+1} := \delta_{a_{n+1}} - \mathbf{br}(\overline{x}_n).\]

\end{itemize}
The associated ODE is the perturbed best response dynamic, given by
\begin{equation}
\label{eq:PBRD}
\dot{x} = \mathbf{br}(x) - x.
\end{equation}
Note that the set of stationary points for this dynamic is exactly PNE, the set of perturbed equilibria. Since the vector field $F$ is taking values in the tangent space  relative to $\Delta$,
$T\Delta := \times T\Delta^i$
the trajectories remain in $\Delta$. By an obvious abuse of language, we will say that a $m \times m$ matrix $A$ is positive definite if, for any $\zeta \in T \Delta$, we have
\[\zeta \neq 0 \Rightarrow \zeta^T A \zeta > 0.\]

In the following, the set of matrices which are positive definite in this sense is denoted $\mathcal{S}^+(T \Delta)$.

\begin{lemma} \label{noisesfp} Assume that for each $i$ the choice function of player $i$ takes values into the interior \footnote{Notice that this property is always satisfied for good stochastic choice functions} of $\Delta^i$.
Then there exists positive values $0< \Lambda^- < \Lambda^+ < + \infty$ and a continuous function $Q : \Delta \rightarrow \mathcal{S}^+(T \Delta) \cap [\Lambda^- I_d, \Lambda^+ I_d]$ such that

\[\mathbb{E} \left(U_{n+1} U_{n+1}^T \mid \mathcal{F}_n \right) = Q (\overline{x}_n).\]
\end{lemma}
\prf
Let, for $x \in \Delta$ and $i \in \{1, \ldots, N\}$ $Q^i(x)$ denote the quadratic form on $T\Delta^i$ defined by
$$Q^i(x)(\zeta^i) = \sum_{\alpha \in A^i} \langle \delta_{\alpha} - br^i(x^{-i}), \zeta^i \rangle^2 br^i(x^{-i})_{\alpha}.$$
Equivalently, $Q^i(x)(\zeta^i)$ is the variance of $\alpha  \mapsto \langle \delta_{\alpha}, \zeta^i \rangle$ under the law $br^i(x^{-i}).$
Let $Q(x)$ denote the quadratic form on $T\Delta$ defined by
$$Q(x)(\zeta) = \sum_{i = 1}^N Q^i(x)(\zeta^i).$$
Since $br^i(x^{-i})_{\alpha} > 0$ and $\{\delta_{\alpha} - br^i(x^{-i}) : \: \alpha \in A^i\}$ spans $T\Delta^i,$ $Q^i(x)$ is non-degenerate for all $i.$ Hence $Q(x)$ is nondegenerate, and by compactness and continuity, there exist $\Lambda^+ \geq \Lambda^- > 0$ such that
$$\Lambda^- \left\|\zeta\right\|^2 \leq Q(x)(\zeta) \leq \Lambda^+ \left\|\zeta\right\|^2, \; \, \forall \zeta \in T \Delta.$$
Now $$\mathbb{E} \left( \langle U_{n+1} U_{n+1}^T \zeta, \zeta \rangle \mid \mathcal{F}_n \right) = Q(\overline{x}_n)(\zeta).$$
Hence the result.
\qed

Finally, the discrete stochastic approximation (\ref{eq:SFP}) is a first case Robbins Monro algorithm with $q=2$, which satisfies hypothesis \ref{noise}.

\subsection{Properties of the best response dynamic for supermodular games}
We assume here that for each $i = 1,\ldots,N$ the action set $A^i$ is equipped with a total ordering denoted $\leq;$ and
we  focus our attention on games such that, for a given player, the reward he obtains by switching to a higher action increases when his opponents choose higher strategies. Such games are called {\em supermodular} and arise in many economic applications : see e.g.  \cite{Top79} or \cite{milRob90}.

\begin{definition} We say that the game $G(N,A,u)$ is (strictly) {\itshape supermodular} if, for any pair of distinct players $(i,j)$ and any action profiles $a = (a^1,..,a^N)$ and $b = (a^1,..,a^N)$ such that $a^i > b^i$ and $a^{-i} = b^{-i}$, the quantity $u^i(a) - u^i(b)$  is (strictly) increasing in $a^j = b^j, \,$ for $j \neq i$.
\end{definition}

\brem   In the particular case where each action set $A^i$ is equal to the couple $\{0,1\}$, the state space is the hypercube $[0,1]^N$ and these games have been defined as {\em coordination games} in (Benaim and Hirch, 1999) \erem

In the remainder of this section we set $A^i = \{1, \ldots, m^i\}$ and we assume that $\leq$ is the natural ordering on integers. For player $i$, we define the invertible linear operator $T^i$:
\[ \Delta^i \rightarrow \mathbb{R}^{m^i - 1}, \; \, (x_k^i)_{k=1,..,m^i} \mapsto ((T^i(x^i))_j)_{j=1,..,m^i-1}\]
with
\[(T^i (x^i))_j = \sum_{k=j+1}^{m^i} x_k^i.\]

\hop Two mixed strategies can be compared via this operator and $T^i(x^i) \leq T^i(y^i)$ if and only if $y^i$ stochastically dominates  $x^i$. In the same spirit, two strategy profiles can be compared introducing the operator $T:$
\[\Delta \rightarrow \times_{i=1,..,N} \mathbb{R}^{m^i-1}, \; \, (x^1,..,x^N) \mapsto (T^1(x^1),..,T^N(x^N)).\]
Naturally, we say that $T(x) \leq T(y)$ if $T^i(x^i) \leq T^i(y^i)$ for $i=1,..,N$ and the order relation relative to $T$ denoted $\leq_T$ in the sequel. The following result is proved in \cite{HofSan02}.

\begin{theorem}[Hofbauer and Sandholm, 2002] \label{HofSan}Assume that the game is strictly supermodular and that every agent plays a smooth fictitious play strategy induced by a good stochastic choice function. Then

\begitem
\iti for $i = 1,..,N, \; \,  y^{-i} \geq_T x^{-i} \, \Rightarrow  \mathbf{br}^i(y^{-i})  \geq_T  \mathbf{br}^i(x^{-i})$.


\itii The conjugate dynamic\footnote{we refer to the dynamic induced by the conjugation relation $T$, defined on $\left\{ (v^1,...,v^N) \in \times_{i=1}^N \mathbb{R}^{m^i-1} \mid \; \, 1 \geq v_1^i \geq ... \geq v^i_{m^i-1} \geq 0 \; \, \forall i \right\}$ and given by $\dot{v} = T \left( br(T^{-1}(v))\right) - v$.} is cooperative and irreducible. Hence, it is strongly monotone. In particular, if $(\mathbf{x}(t))_{t \geq 0}$ and $(\mathbf{y}(t))_{t \geq 0}$ solve (\ref{eq:PBRD}) with $\mathbf{x}(0) \leq_T \mathbf{y}(0)$ (and $\mathbf{x}(0) \neq \mathbf{y}(0)$) then, for any $t \geq 0$, $\mathbf{x}(t) \leq_T \mathbf{y}(t)$,

\itiii There exists two perturbed equilibria $\underline{x} \leq_T \overline{x}$ such that any chain recurrent set relative to the perturbed best response dynamic is included into the interval  $[\underline{x},\overline{x}]$,

\end{itemize}
\end{theorem}
\zdeux
Hofbauer and Sandholm then used this theorem combined with results from \cite{Ben00} to describe the limit set  of  stochastic fictious plays for supermodular game.
In view of the new results obtained in this paper and specifically in section \ref{sec:coop} we are are now able to improve notably their results and to prove the convergence of stochastic fictious play for supermodular games in full generality.
\begin{theorem} Assume that the assumptions of previous theorem are satisfied. Then the limit set of  $(\overline{x}_n)_n$ is almost surely an ordered arc of PNE that are not linearly unstable. If we furthermore assume that the choice function is real analytic (for instance in the logit case), then  $(\overline{x}_n)_n$ almost surely converges toward a non linearly unstable PNE.
\end{theorem}
\prf By Lemma \ref{noisesfp} and Theorem \ref{HofSan}, the conditions to apply Corollary \ref{th:cooprm} are met. \qed

\section{Appendix}
\subsection{Proof of Proposition \ref{th:nonconvsde}}
The assumptions on $\gamma$ easily imply that $$\frac{\gamma(t)}{\gamma(s+t)} \geq \frac{1}{\gamma(s)} \geq e^{bs}.$$ Thus
 $$\omega(t, a \sqrt{\gamma(t)}, T) \leq C\int_0^{\infty} \exp{(- a^2 e^{bs} C(T))}$$ and condition $(iii)$ of hypothesis \ref{H2}holds.
Let $$A_s^t = [DF(\Phi_s(X_t)) - \frac{1}{2} \frac{\dot{\gamma}(t+s)}{\gamma(t+s)}]$$ and
let  $\{Y_s^t, s \geq 0\}$ be solution to
$$dY_s^t = A_s^t  Y_s^t + dB_{t+s}$$
with initial condition $Y_0^t = 0.$
Condition $(i)$ of Hypothesis \ref{H2} follows from the following lemma.
\blem \label{diff (i)}
$$\lim_{t \rar \infty} \Pr(\sup_{0 \leq s \leq 1} \|Y_s^t - \frac{X_{t+s} - \Phi_s(X_t)}{\sqrt{\gamma(t+s)}}\| \geq \eps |\cF_t) = 0.$$ In particular, Hypothesis \ref{H2} $(i)$ holds with $Y(t) = Y^{t-1}_1$ for all $t \geq 1.$
\elem
\prf
Set $\alpha(s) = 1/\sqrt{\gamma(s)}, \,  Z^t_s = X_{t+s} - \Phi_s(X_t)$ and  $\hat{Y}_s^t = \alpha(t+s) Z_s^t$. Then
$$dZ_s^t = (F(X_{t+s}) - F(\Phi_s(X_t))) ds + \sqrt{\gamma(t+s)} dB_{t+s}$$
$$= [DF(\Phi_s(X_t)) Z_s^t  + o(\|Z_s^t\|)]ds + \sqrt{\gamma(t+s)} dB_{t+s}.$$
Hence
$$d\hat{Y}_s^t = [DF(\Phi_s(X_t)) + \frac{\dot{\alpha}(t+s)}{\alpha(t+s)}] \hat{Y}_s^t + dB_{t+s} + \alpha(t+s) o(\|Z_s^t\|),$$
where $o(z) = z \eta(z)$ and $\lim_{z \rar 0} \eta(z) = \eta(0) = 0.$
Then
$$Y_s^t - \hat{Y}_s^t = \int_0^s A_u^t(Y_u^t - \hat{Y}_u^t) du + \int_0^s \alpha(t+u)o(\|Z_u^t\|)du.$$
Thus, by Gronwall's inequality,
$$\sup_{0 \leq s \leq 1} \|Y_s^t - \hat{Y}_s^t\| \leq e^K R_t $$
with
$$R_t = \sup_{0 \leq s \leq 1} \alpha(t+s) o(\|Z_s^t\|)$$ and
\beq
\label{defK}
K = \sup_{s,t} \|A_s^t\| \leq \|DF\| + \frac{a}{2}.
\eeq To conclude
the proof  it remains to show that
 $$\Pr(R_t \geq \delta |\cF_t) \rar 0$$ as $t \rar \infty.$

It follows from the estimate given in example \ref{ex:diffusion} that $$\Pr(\sup_{0 \leq s \leq 1} \|Z^t_s\| \geq \delta | \cF_t) \leq \int_t^{t+1}C \exp (\frac{-\delta^2 C(1)}{\gamma(s)})ds \leq C \exp(-\frac{\delta^2 C(1)}{\gamma(t+1)})$$
Thus
\begin{eqnarray*}
\Pr(\sup_{0 \leq s \leq 1} \alpha(t+s) \|Z_s^t\|\geq R | \cF_t) &\leq&  \Pr(\|Z_s^t\|\geq \frac{R}{\alpha(t+1)}| \cF_t))\\
&\leq&  C \exp (- R^2 C(1)).
\end{eqnarray*}
Now,
$$\Pr(\sup_{0 \leq s \leq 1} \alpha(t+s)\|Z_s^t\| \eta(\|Z_s^t\|) \geq \delta |\cF_t)$$
$$ \leq
\Pr(\sup_{0 \leq s \leq 1} \alpha(t+s)\|Z_s^t\| \geq R |\cF_t) + \Pr(\sup_{0 \leq s \leq 1} \eta(\|Z_s^t\|) \geq \frac{\delta}{R} |\cF_t).$$
$$\leq C \exp (- R^2 C(1)) + \Pr(\sup_{0 \leq s \leq 1} \eta(\|Z_s^t\|) \geq \frac{\delta}{R}|\cF_t).$$
Since $\lim_{z \rar 0} \eta(z) = 0,$
$$\limsup_{t \rar \infty} \Pr(\sup_{0 \leq s \leq 1} \alpha(t+s)\|Z_s^t\| \eta(\|Z_s^t\|) \geq \delta |\cF_t) \leq C \exp (- R^2 C(1))$$ and since $R$ is arbitrary, this proves the result.
\qed

It remains to prove that condition $(ii)$ of hypothesis \ref{H2} holds.
\blem
\label{easygauss}
Let $\Sigma$ be a $n \times n$ self-adjoint positive definite matrix and $$f_{\Sigma}(x) = \frac{\exp{(-\frac{1}{2}\langle \Sigma^{-1} x, x \rangle})}{\sqrt{\det(\Sigma) (2\pi)^n}}$$ the density of a centered Gaussian vector with covariance $\Sigma.$ Let $0 < \alpha
\leq \beta$ respectively denote the smallest and largest eigenvalues of $\Sigma.$  Then
$$f_{\Sigma}(x) \geq (\frac{\alpha}{\beta})^{n/2} f_{\alpha Id}(x).$$
\elem
\prf Follows from the estimates $\det(\Sigma) \leq \beta^n$ and
$\langle \Sigma^{-1} x, x \rangle \leq\frac{ \|x\|^2}{\alpha}.$ \qed

Since $Y^t_s$ is a linear function of $\{B_{t+u},  \; 0 \leq u \leq s\}$, it is a Gaussian vector under the conditional probability $\Pr( \cdot | \cF_t)$. By Ito's formulae, its covariance matrix is solution to
$$\frac{d\Sigma^t_s}{ds} = A_s^t \Sigma^t_s + \Sigma^t_s A_s^{t*} + Id$$ with initial condition $\Sigma^t_0 = 0;$ where $A_s^{t*}$ stands for the transpose of $A_s^t.$ It is then easy to check that
$$\Sigma^t_s = \int_0^s U^t(u) U^{t*}(u) du$$ where $U^t(s)$ is the solution to
\beq
\label{odeU}
\frac{dU}{ds} = A_s^t U, \,  U(0) = Id.
\eeq
Using (\ref{odeU}) we see that  $U^t(s)$ is invertible and that its inverse $(U^t(s))^{-1}$ solves
$$\frac{dV}{ds} = -V A_s^t, \, V(0) = Id.$$
Using again (\ref{odeU}) combined with the estimate (\ref{defK}) and Gronwall's lemma, we get
$$\|U^t(s)\| \leq e^{Ks}.$$ Similarly
$$\|(U^t(s))^{-1}\| \leq e^{Ks}.$$  It follows that
for all vector $h,$ $$e^{-Ks} \|h\| \leq \|U^t(s)h\| \leq e^{Ks}\|h\|.$$
Hence
$$a \|h\|^2 \leq \langle \Sigma^t_1 h, h \rangle \leq  b \|h\|^2,$$
where $a = \int_0^1 e^{-2Ku} du$ and  $b = \int_0^1 e^{2Ku} du.$
The result then follows from Lemma (\ref{easygauss}). $\; \; \blacksquare$

\subsection{Proof of Proposition \ref{RMNC}}
Recall that $(\mathcal{F}_n)_n$ is a given filtration to which the stochastic process $(x_n)_n$ is adapted. Let $m_n := \sup \{k \in \mathbb{N} \mid \tau_k \leq n \}$ and call $(\mathcal{G}_n)_n$ the sigma algebra $(\mathcal{F}_{m_n})_n$. Let $n \geq 1$ and $k_n := m_{n+1} - m_n$.  We denote by $t^n_j$ the quantity $\tau_{m_n + j} - \tau_{m_n}$ ($j=0,..,k_n$) and $t_n := t^n_{k_n}$. Notice that $|t_n -1| \leq \gamma_{m_n}$.

\hop For the continuous time interpolated process induced by a discrete process $(x_n)_n$, hypothesis \ref{H2} is satisfied if there exists a vanishing positive sequence $(\gamma(n))_n$ and a $\mathcal{G}_n$-adapted random sequence $(Y_n)_n$ such that

\begitem
\item[$(i)$] for any $\alpha >0$,
\[\lim_{n \rightarrow + \infty} \, \mathbb{P} \left(\left\|\frac{x_{m_{n+1}} - \Phi_{t_n} (x_{m_n})}{\sqrt{\gamma(n)}} - Y_{n+1} \right\| > \alpha \mid \mathcal{G}_n  \right) = 0,\]

\item[$(ii)$] for any open set $O \subset \mathbb{R}^d$, there exists a positive number $\delta$ such that
\[\liminf_{n \rightarrow + \infty} \, \mathbb{P} \left(Y_{n+1} \in O \mid \mathcal{G}_n \right) > \delta \; \, \mbox{ almost surely.}\]

\item[$(iii)$] there exists $a>0$ such that
\[\limsup_{n \rightarrow + \infty} \,  \omega (n,a\sqrt{\gamma(n)},T) <1.\]
\end{itemize}

\hop  Let $\gamma(n) := \sum_{k=1}^{k_n} \gamma_{m_n + k}^2$. First, by proposition \ref{omega}, the map $\omega$ corresponding to  the process $(x_n)_n$ is given by

\[\omega (n,\delta,T) = \frac{B \int_n^{+ \infty} \overline{\gamma}(u) du}{\delta^2}.\]

\hop Hence,

\[\omega (n,a \sqrt{\gamma(n)},T) \leq \frac{B}{a^2} \frac{\sum_{m_n}^{+ \infty} \gamma_i^2}{\sum_{m_n+1}^{m_{n+1}} \gamma_i^2}.\]
Since
\[\limsup_n \frac{\sum_{m_n}^{+ \infty} \gamma_i^2}{\sum_{m_n+1}^{m_{n+1}} \gamma_i^2}  < + \infty,\]
the quantity $\omega (n,a \sqrt{\gamma(n)},T)$ is smaller than $1$, for $a$  large enough. The next lemma corresponds to Lemma  \ref{diff (i)}.

\begin{lemma}  \label{RM (i)}
Point $(i)$ is satisfied for this choice of $(\gamma(n))_n$ and the random sequence $(Y_n)_n$ given by
\[\frac{1}{\sqrt{\gamma(n-1)}} \sum_{j=1}^{k_{n-1}} \gamma_{m_{n-1}+j} \left(\prod_{k=j+1}^{k_{n-1}} \left(I_d + \gamma_{m_{n-1} + k} DF (\phi_{t^{n-1}_{k-1}} (x_{m_{n-1}}))\right)\right) U_{m_{n-1} +j}.\]
\end{lemma}
\zun

\hop {\bfseries Proof.} Set $\hat{Y}_{n+1} := \frac{x_{m_{n+1}} - \phi_{t_n} (x_{m_n})}{\sqrt{\gamma(n)}}$.  We have, for $j=0,..,k_n -1$,
\[\phi_{t^n_{j+1}}(x_{m_n}) - \phi_{t^n_j}(x_{m_n}) = \gamma_{m_n + j+1} F \left(\phi_{t^n_j} (x_{m_n}) \right) + \mathcal{O} (\gamma_{m_n + j}^2).\]

\hop Then, denoting
\[\hat{Y}^n_j := \frac{1}{\sqrt{\gamma(n)}} \left( x_{m_n + j} - \phi_{t^n_j} (x_{m_n}) \right) (j=0,..,k_n),\]
we have
\begin{eqnarray*}
\hat{Y}^n_{j+1} - \hat{Y}^n_j &=& \frac{\gamma_{m_n + j + 1}}{\sqrt{\gamma(n)}} \left[ F (x_{m_n + j})  - F \left(\phi_{t^n_j} (x_{m_n}) \right)  +  U_{m_n + j + 1} \right] \\
&&+ \mathcal{O}  \left( \frac{\gamma_{m_n + j +1}^2}{\sqrt{\gamma(n)}}\right).
\end{eqnarray*}
Consequently,
\begin{eqnarray*}
\hat{Y}^n_{j+1} - \hat{Y}^n_{j} &=& \gamma_{m_n +j+1} \left( DF \left(\phi_{t^n_j}(x_{m_n}) \right)  \hat{Y}^n_{j} + \frac{R^n(j)}{\sqrt{\gamma(n)}} + \frac{U_{m_n+j+1}}{\sqrt{\gamma(n)}} \right) \\
&& + \mathcal{O}  \left( \frac{\gamma_{m_n + j +1}^2}{\sqrt{\gamma(n)}}\right)\; \; j = 0,..,k_n-1,
\end{eqnarray*}

\hop where
\[R^n(j) := F(x_{m_n +j}) - F(\phi_{t^n_j} (x_{m_n})) - DF \left(\phi_{t^n_j}(x_{m_n}) \right) \cdot \left(x_{m_n + j} - \phi_{t^n_j} (x_{m_n}) \right).\]
By a recursive argument,
\begin{eqnarray*}
&& \hat{Y}_{n+1} -Y_{n+1} = \hat{Y}^n_{k_n} - Y_{n+1} =\\
&& \frac{1}{\sqrt{\gamma(n)}} \sum_{j=1}^{k_n} \gamma_{m_n+j} \left(\prod_{k=j+1}^{k_n} \left(I_d + \gamma_{m_n + k} DF (\phi_{t^n_{k-1}} (x_{m_n}))\right)\right) R^n(j)\\
&& + \mathcal{O}\left(e^{-n/2}\right).
\end{eqnarray*}

\hop since $\hat{Y}^n_0 = 0$ and $\sum_{j=0}^{k_n-1} \frac{\gamma_{m_n + j + 1}}{\sqrt{\gamma(n)}} = \sqrt{\gamma(n)} = \mathcal{O}(e^{-n/2})$.

Recall that $\sum_{j=1}^{k_n} \gamma_{m_n+j} \leq 1 + \gamma_{m_{n+1}}$ and $DF$  is bounded.   Consequently, there exists a real number  $K$ such that for $n$ large enough,
\begin{eqnarray*}
&&\frac{1}{\sqrt{\gamma(n)}} \left\| \sum_{j=1}^{k_n} \gamma_{m_n+j} \left(\prod_{k=j+1}^{k_n} \left(I_d + \gamma_{m_n + k} DF (\phi_{t^n_{k-1}}(x_{m_n}))\right)\right) R^n(j)\right\| \\
&& \leq e^K \frac{1}{\sqrt{\gamma(n)}} \sup_{j=1,..,k_n} R^n(j) = e^K R_n,
\end{eqnarray*}
where  $R_n := \frac{1}{\sqrt{\gamma(n)}} \sup_{j=1,..,k_n} R^n(j)$.
By an application of results due to Benaim (see \cite{Ben99}, proposition 4.1, formula (11) and identity (13) with $q=2$), we have
\[\mathbb{E} \left( \sup_{j=0,..,k_n-1} \|x_{m_n+j} - \phi_{t^n_j} (x_{m_n})\|^2 \mid \mathcal{G}_n \right) \leq C \gamma(n).\]

\hop where $C$ is some positive constant. Additionally, by definition of $DF$,
\[R^n(j)^2 \leq h \left(\|x_{m_n+j} - \phi_{t^n_j} (x_{m_n}) \|^2 \right),\]

\hop for some function $h: \mathbb{R}_+^* \rightarrow \mathbb{R}_+^*$, strictly increasing and such that $h(x)/x \rightarrow_{x \rightarrow 0^+} 0^+$. An immediate consequence is that
\begin{eqnarray*}
&& \mathbb{P} \left( R_n \geq \alpha \mid \mathcal{G}_n \right) \\
&\leq& \mathbb{P} \left( \sup_{j=0,..,k_n-1} h \left(\|x_{m_n+j} - \phi_{t^n_j} (x_{m_n})\|^2\right) \geq \alpha^2 \gamma(n) \mid \mathcal{G}_n \right)\\
&\leq& \mathbb{P} \left(\sup_{j=0,..,k_n -1}\|x_{m_n+j} - \phi_{t^n_j} (x_{m_n})\|^2 \geq h^{-1} \left(\alpha^2 \gamma(n) \right) \mid \mathcal{G}_n\right)\\
&\leq& \frac{C \gamma(n)}{h^{-1} \left( \alpha^2 \gamma(n) \right)} \rightarrow_{n \rightarrow + \infty} 0,
\end{eqnarray*}

\hop which proves the result. $\; \; \blacksquare$

\zdeux

\hop To simplify notations, we call $E$ the euclidian space $\mathbb{R}^d$. Given $n \in \mathbb{N}$, the random variable $x_n$ can be written $h_n(U_1,..,U_n)$, where $h_n : (E^{n}, (\mathcal{B}_E)^{n}) \rightarrow (E, \mathcal{B}_E)$ is a measurable function. We denote by $\mathcal{P}_U$ the probability distribution induced by the measurable process $U=(U_n)_n : (\Omega, \mathcal{F}) \rightarrow (E^{\NN}, (\mathcal{B}_E)^{\NN})$. We keep the notation  $\mathcal{F}_n$ for the sigma field $(\mathcal{B}_E)^n \times E^{\NN}$ when it does not imply any ambiguity.

\begin{proposition}  There exists a function $P_n: (\mathcal{B}_{E})^{\NN} \times E^{\NN} \rightarrow [0,1]$ called a \emph{regular conditional distribution of  $U$ given $\mathcal{F}_n$} in the sense that, for any $u \in E^{\NN}$, $P_n(\cdot,u)$ is a probability measure on $((\mathbb{R}^d)^{\NN}, (\mathcal{B}_E)^{\NN})$ and that, for any $B \in (\mathcal{B}_E)^{\NN})$, the random variable $P_n(B,\cdot)$ is $\mathcal{F}_n$-measurable with
\[\mathbb{P}_n(B,\cdot) = \mathbb{P}_U(B \mid \mathcal{F}_n)(\cdot) \; \, \mathbb{P}_U-\mbox{almost surely.}\]
\end{proposition}

For convenience, given $u \in E^{\NN}$, we denote by $\mathbb{P}_n^{u}$ the probability measure $\mathbb{P}_n(\cdot,u)$ and $\mathbb{E}_n^{u}$ the corresponding expectation. Given a measurable function $y : (E^{\NN}, (\mathcal{B}_{E})^{\NN}) \rightarrow (E,\mathcal{B}_E)$, we have
\[\mathbb{E}^{\omega}_n(y) = \mathbb{E}_U(y\mid \mathcal{F}_n) = \mathbb{E}(y(U) \mid \mathcal{F}_n) \; \, \mathbb{P}_U-\mbox{a.s.}.\]

\begin{lemma} Let $k < i$ be two natural numbers and $y : (E^{\NN}, (\mathcal{B}_{E})^{\NN}) \rightarrow (E,\mathcal{B}_E)$ be a measurable function. There exists a subset $\Omega_0(y) \subset E^{\NN}$ such that $\mathbb{P}_U(\Omega_0(y)) = 1$ and, for any $u_0 \in \Omega_0(y)$, $\mathbb{E}^{u_0}_k (y\mid \mathcal{F}_i)$ and $\mathbb{E}_U(y\mid \mathcal{F}_i)$ are $\mathbb{P}_U$-almost surely equal.
\end{lemma}

\hop {\bfseries Proof.} The random variable $z := \mathbb{E}_U(y \mid \mathcal{F}_i)$ is $\mathcal{F}_i$-measurable. Pick a countable $\pi$-class $\mathcal{D}$ such that $\sigma(\mathcal{D}) = \mathcal{F}_k$. Given $A \in \mathcal{D}$, there exists a set $\Omega_0(y,A)$ such that $\mathbb{P}_U(\Omega_0(y,A))=1$ and, for any $u_0 \in \Omega_0(y,A)$, we have
\begin{itemize}
\item[$(1)$] $\mathbb{E}^{u_0}_k \left( \mathbb{E}(\mathbb{I}_A y \mid \mathcal{F}_i)\right) = \mathbb{E}_U \left(\mathbb{E}(\mathbb{I}_A y \mid \mathcal{F}_i) \mid \mathcal{F}_k \right)(u_0),$
\item[$(2)$] $\mathbb{E}^{u_0}_k \left( \mathbb{I}_A y\right) = \mathbb{E}_U \left(\mathbb{I}_A y  \mid \mathcal{F}_k \right)(u_0).$
\item[$(3)$] $ \mathbb{I}_A \mathbb{E}_U(y \mid \mathcal{F}_i) =  \mathbb{E}_U(\mathbb{I}_A y \mid \mathcal{F}_i)$ $\; \mathbb{P}_k^{u_0}$-a.s.
\end{itemize}

Let us construct $\Omega_0(y,A)$. First, there exist two sets $\Omega_0^1(y,A)$ and $\Omega_0^2(y,A)$ on which respectively points $(1)$ and $(2)$ are satisfied and such that $\mathbb{P}_U(\Omega_0^j(y,A)) = 1, \, j=1,2$.
Now for the last point, one must first consider a set $\Omega^3(y,A)$ such that $\mathbb{P}_U(\Omega^3(y,A))=1$ and, for any $u \in \Omega^3(y,A)$,
\[\mathbb{I}_A(u) \mathbb{E}_U(y \mid \mathcal{F}_i)(u) =  \mathbb{E}_U(\mathbb{I}_A y \mid \mathcal{F}_i)(u)\]
Then, by definition of $\mathbb{P}^{u_0}_k$, there exists a set  $\Omega_0^3(y,A)$ (which depends on $\Omega^3(y,A)$) such that, $\mathbb{P}_U(\Omega_0^3(y,A))=1$ and, for any $u_0 \in \Omega_0^3(y,A)$,
\[\mathbb{P}^{u_0}_k(\Omega^3(y,A)) = \mathbb{P}_U(\Omega^3(y,A) \mid \mathcal{F}_k)(u_0) = 1. \]
Finally, pick $\Omega_0(y,A) := \Omega^1_0(y,A) \cap \Omega^2_0(y,A) \cap \Omega^3_0(y,A)$.

Now take
\[\Omega_0(y) := \cap_{A \in \mathcal{D}} \Omega(y,A).\]
By countability of $\mathcal{D}$, we have $\mathbb{P}_U(\Omega_0(y)) = 1$. There remains to prove that, for any $u_0 \in \Omega_0(y)$,
\[\int_A z \, d\mathbb{P}^{u_0}_k = \int_A y \, d\mathbb{P}^{u_0}_k, \; \, \mbox{for any} \,  A \in \mathcal{D}.\]

\begin{eqnarray*}
\mathbb{E}_k^{u_0}(\mathbb{I}_A z) &=& \mathbb{E}^{u_0}_k \left( \mathbb{I}_A \mathbb{E}_U(y \mid \mathcal{F}_i)\right)\\
&=& \mathbb{E}^{u_0}_k \left(\mathbb{E}_U(\mathbb{I}_A y \mid \mathcal{F}_i) \right) \\
&=& \mathbb{E}_U \left(\mathbb{E}_U(\mathbb{I}_A y \mid \mathcal{F}_i) \mid \mathcal{F}_k \right)(u_0) \\
&=& \mathbb{E}_U \left( \mathbb{I}_A y \mid \mathcal{F}_k \right)(u_0)\\
&=& \mathbb{E}^{u_0}_k \left( \mathbb{I}_A y \right).
\end{eqnarray*}
The second equality follows from point $(3)$, the third from point $(1)$ and the fifth from point $(2)$. The lemma is proved. $\; \; \blacksquare$
\zdeux

The following result is due to \cite{HalHey80} (see Theorem 3.4 or Theorem 2 page 351 in \cite{ChoTei98} for a version adapted to our situation). It is a central limit result for double arrays. We apply it to prove point $(ii)$.

\begin{theorem}[Hall and Heyde] \label{th:HH} For any $n \geq 1$, let $k_n$ be a positive integer and $(\Omega_n, \mathcal{F}^n, \mathbb{P}_n)$ a probability space. Consider $\mathcal{F}^n_1 \subset \mathcal{F}^n_2 \subset ... \subset \mathcal{F}^n_{k_n} \subset \mathcal{F}^n$ an increasing family of sigma fields and  $(y^n_j)_{j=1,..,k_n}$ a $(\mathcal{F}^n_j)_{j=1,..,k_n}$-adapted family of random variables. Assume that

\begitem

\item[$*$] for $j=1,..,k_n$,
\[\mathbb{E}_n \left(y_j^n \mid \mathcal{F}^n_{j-1} \right) = 0,\]

\item[$*$] we have
\[\sum_{j=1}^{k_n} \mathbb{E}_n \left(\|Y^n_j\|^2 \mathbb{I}_{\|Y^n_j\| > \varepsilon} \mid \mathcal{F}^n_{j-1}\right) \xrightarrow[n \rightarrow + \infty]{dist.} 0,\]

\item[$*$] there exists a positive, $\mathcal{F}^n_1$-adapted random sequence $(w_n)_n$ such that
\[\sum_{i=1}^{k_n} \mathbb{E}_n \left(y^n_j  (y^n_j)^T \mid \mathcal{F}^n_{j-1} \right) - w_n \xrightarrow[n \rightarrow + \infty]{dist.} 0,\]

\item[$*$] there exists a positive random matrix $\eta$, defined on some probability space $(\Omega, \mathcal{F},\mathbb{P})$, which satisfies
\[\sum_{j=1}^{k_n} \mathbb{E}_n \left(y^n_j  (y^n_j)^T \mid \mathcal{F}^n_{j-1} \right)  \xrightarrow[n \rightarrow + \infty]{dist.} \eta.\]

\end{itemize}

\hop Then, denoting $y_{n+1} := \sum_{j=1}^{k_n} y^n_j$, the sequence $(y_n)_n$ converges in distribution to some random variable $y$ defined on $(\Omega,\mathcal{F},\mathbb{P})$ and whose characteristic function is given by $\mathbb{E} \left(e^{-\frac{1}{2} t^T \eta t}\right)$. In particular,
\[ \lim_{n \rightarrow + \infty}\mathbb{E}_n \left( e^{i<t,y_{n+1}>} \right) = \mathbb{E} \left( e^{-\frac{1}{2} t^T \eta t}\right).\]

\end{theorem}
\ztrois

Let us get back to our settings. Let $n \in \mathbb{N}$ and $j \in \{1,..k_n\}$. Consider the measurable functions $y^n_j : (E^{\NN}, (\mathcal{B}_{E})^{\NN}) \rightarrow (E,\mathcal{B}_E)$, given by
\[y^n_j(u) := \frac{\gamma_{m_n+j} }{\sqrt{\gamma(n)}}  \left(\prod_{k=j+1}^{k_n} \left(I_d + \gamma_{m_n + k} DF (\phi_{t^n_{k-1}} (x_{m_n}))\right)\right) u_{m_n +j},\]
where $x_n = h_n(u_1,..,u_n)$. Finally, call $y_n := \sum_{j=1}^{k_n} y^n_j$

\begin{corollary} Given  a nonempty open set $O$ in $E$, there exist $\delta > 0$ and a set $\Omega_0$ such that $\mathbb{P}_U(\Omega_0) =1$ and, for any $u_0 \in \Omega_0$,
\[\liminf_n \mathbb{P}_U \left(y_{n+1} \in O \mid \mathcal{G}_n \right)(u_0)  > \delta. \]
\end{corollary}
\zdeux

\hop {\bfseries Proof.}  Let $\Omega_0$ be the set
\[\bigcap_{n \in \mathbb{N},j=1,..,k_n, r \in \mathbb{Q}} \Omega_0\left(y^n_j, \, \|y_j^n\|^2 \mathbb{I}_{\|y^n_j\| > r}), \,  y^n_j (y_j^n)^T, \, \mathbb{I}_{\|x_{m_n+j}- \Phi_{t^n_j}(x_{m_n})\|>r}, \mathbb{I}_{y_n \in O}  \right).\]
By countability,  $\mathbb{P}(\Omega_0)=1$. Pick $u_0 \in \Omega_0$. We apply Theorem \ref{th:HH} to  $(\Omega_n,\mathcal{F}_n,\mathbb{P}_n) := (E^{\NN},(\mathcal{B}_E)^{\NN},\mathbb{P}^{u_0}_{m_n})$,
$\mathcal{F}^n_j = \mathcal{F}_{m_n+j}$ and the double array of random variables $(y_j^n)_{n,j}$.
\zdeux

We now verify that the assumptions required to apply Theorem \ref{th:HH} hold. First of all
\[\mathbb{E}_{m_n}^{u_0} \left(y^n_j \mid \mathcal{F}^n_{j-1}\right) = \mathbb{E}_U(y_j^n \mid \mathcal{F}^n_{j-1}) = 0 \; \, \mbox{a.s}.\]
\zun

\hop Secondly, let
\[\Pi_{n,j} := \prod_{k=j+1}^{k_n} \left(I_d + \gamma_{m_n + k} DF \left( \phi_{t^n_{k-1}} (x_{m_n})\right) \right).\]
A simple computation gives
\[ e^{-2 \|DF\|_{\infty}}\leq \left\|\Pi_{n,j} \right\|\leq e^{\|DF\|_{\infty}}.\]

Recall that there exists $p>1$ such that the sequence of random variables $(\mathbb{E}_U\left(\|u_n\|^{2p} \mid \mathcal{F}_{n-1}\right))_n$ is almost surely bounded. Hence, taking $q$ such that $1/p + 1/q =1$ and choosing $r \in \mathbb{Q}$,
\begin{eqnarray*}
\mathbb{E}_U \left(\|y^n_j\|^2 \mathbb{I}_{\|y^n_j\| > r} \mid \mathcal{F}^n_{j-1}\right) &\leq& \mathbb{E}_U\left(\|y^n_j\|^{2p} \mid \mathcal{F}^n_{j-1}\right)^{1/p} \mathbb{P}_U \left(\|y^n_j\|^{2p} > r^{2p} \mid \mathcal{F}^n_{j-1}\right)^{1/q} \\
&\leq& \frac{1}{r^{2p/q}}\mathbb{E}_U\left(\|y^n_j\|^{2p} \mid \mathcal{F}^n_{j-1}\right)\\
&\leq& \frac{1}{r^{2p/q}} \frac{\gamma_{m_n+j}^{2p}}{\gamma(n)} e^{\|DF\|_{\infty}} \mathbb{E}_U\left(\|u_{m_n+j}\|^{2p} \mid \mathcal{F}^n_{j-1}\right)\\
&\leq& C(r) \frac{\gamma_{m_n+j}^{2p}}{\gamma(n)} \mathbb{E}_U\left(\|u_{m_n+j}\|^{2p} \mid \mathcal{F}^n_{j-1}\right).
\end{eqnarray*}

Consequently,
\[\sum_{j=1}^{k_n} \mathbb{E}_U \left(\|y^n_j\|^2 \mathbb{I}_{\|y^n_j\| > r} \mid \mathcal{F}^n_{j-1}\right) \leq C(r) \sup_j \gamma_{m_n+j}^{2(p-1)} \sup_j \mathbb{E}_U\left(\|u_{m_n+j}\|^{2p} \mid \mathcal{F}^n_{j-1}\right), \]
which converges to $0$ almost surely.
 Since $u_0$ belongs to the set $\Omega_0(\|y_j^n\|^2 \mathbb{I}_{\|y^n_j\|>r})$, for any $j=1,..k_n$,
 \[\sum_{j=1}^{k_n} \mathbb{E}^{u_0}_{m_n}\left(\|y^n_j\|^2 \mathbb{I}_{\|y^n_j\| > r}\mid \mathcal{F}^n_{j-1} \right) = \sum_{j=1}^{k_n} \mathbb{E}_U \left(\|y^n_j\|^2 \mathbb{I}_{\|y^n_j\| > r} \mid \mathcal{F}^n_{j-1}\right) \; \, \mathbb{P}_U-\mbox{a.s.}\]
 and the second point holds.
\zun

\hop  From now on, we call
\[W_n := \sum_{j=1}^{k_n } \mathbb{E}_U \left( \left(y^n_j \right)  \left(y^n_j  \right)^T \mid \mathcal{F}_{n,j-1} \right).\]

\hop We have
\begin{eqnarray*}
&&\mathbb{E}_U \left( \left(y^n_j \right)  \left(y^n_j \right)^T \mid \mathcal{F}_{n,j-1} \right) \\
&=& \frac{1}{\gamma(n)} \gamma^2_{m_n + j}  \Pi_{n,j} \mathbb{E}_U \left( u_{m_n + j} u_{m_n + j}^T \mid \mathcal{F}_{n,j-1}\right)  \Pi_{n,j}^T\\
&=& \frac{1}{\gamma(n)} \gamma_{m_n + j}^2 \Pi_{n,j} Q(x_{m_n+j-1}) \Pi_{n,j}^T.
\end{eqnarray*}

\hop Consequently,
\[W_n = \frac{1}{\gamma(n)} \sum_{j=1}^{k_n} \gamma_{m_n + j}^2 \Pi_{n,j} Q(x_{m_n+j-1}) \Pi_{n,j}^T.\]

\hop Let $w_n$ be the $\mathcal{F}_{n,1}$-measurable random variable  defined by
\[w_n := \frac{1}{\gamma(n)} \sum_{j=1}^{k_n}  \gamma_{m_n + j}^2 \Pi_{n,j} Q\left(\phi_{t^n_{j-1}} (x_{m_n})\right) \Pi_{n,j}^T.\]

Pick $r \in \mathbb{Q}$. By definition of $\Omega_0$ and assumption \ref{H1} $(i)$,
\begin{eqnarray*}
&&\mathbb{P}_{m_n}^{u_0} \left(\sup_{j=1,..k_n} \left\|\phi_{t_{j-1}^n}(x_{m_n}) - x_{m_n+j-1} \right\| >r  \right) \\
&& = \mathbb{P} \left(\sup_{j=1,..k_n} \left\|\phi_{t_{j-1}^n}(x_{m_n}) - x_{m_n+j-1} \right\| >r  \mid \, \mathcal{G}_n \right) \\
&& \leq \omega(n,r,1) \rightarrow 0,
\end{eqnarray*}
which implies that
\[W_n - w_n \xrightarrow[n \rightarrow + \infty]{dist} 0\]

Since the application $Q$ takes values in $[\Lambda^- I_d, \Lambda^+ I_d]$ and $\|\Pi_{n,j}\|$ is bounded above and away from zero, we have
\[0 < a^- \leq \Pi_{n,j} Q(x_{m_n+j-1}) \Pi_{n,j}^T \leq a^+ < + \infty.\]

\hop $W_n$ is a convex combination of such quantities, therefore is bounded. Pick some increasing sequence of integers $(n_k)_k$.
$(W_{n_k})_k$ admits a subsequence $(W_{n'_k})_k$ which converges in distribution to some random variable $\eta^{u_0}$, defined on the probability space induced by $U$ and which takes values in $\mathcal{S}^+(\mathbb{R}^d) \cap [a^- I_d, a^+ I_d]$.

\hop Now by Theorem \ref{th:HH} ,
\[y_{n'_k}  \xrightarrow[n \rightarrow + \infty]{\mathcal{L}} y^{u_0},\]
\hop with $\mathbb{E}_U(e^{i<t,y^{u_0}>}) = \mathbb{E} \left( e^{- \frac{1}{2} t^T \eta^{u_0} t}\right)$.
In particular, by definition of $\Omega_0$,
\[\lim_k \mathbb{P}_U \left(y_{n'_k+1} \in O \mid \mathcal{G}_{n'_k}) \right)(u_0) = \lim_k \mathbb{P}^{u_0}_{n'_k}(y_{n'_k+1} \in O)  = \mathbb{P}(y^{u_0} \in O) > \delta,\]
where $\delta$ depends on the parameters $a$ and $b$ but not on $u_0 \in \Omega_0$ and $(n'_k)_k$. The proof is complete. $\; \; \blacksquare$
\zun

\bibliographystyle{apalike}

\bibliography{supermodular}

\end{document}